\pdfoutput=1
\RequirePackage{ifpdf}
\ifpdf 
\documentclass[pdftex]{sigma}
\else
\documentclass{sigma}
\fi

\usepackage{enumitem}

\usepackage[mathscr]{eucal}

\usepackage{dsfont}

\usepackage[all,cmtip]{xy}
\usepackage{tikz-cd} 

\numberwithin{equation}{section}

\newtheorem{Theorem}{Theorem}[section]
\newtheorem*{Theorem*}{Theorem}

\newtheorem{Proposition}[Theorem]{Proposition}

\theoremstyle{definition}
\newtheorem{Definition}[Theorem]{Definition}
\newtheorem{Example}[Theorem]{Example}
\newtheorem{Remark}[Theorem]{Remark}

\newcommand{\mto}{\mapsto}
\newcommand{\be}{\begin{array}}
\newcommand{\ee}{\end{array}}

 \newcommand\vf\varphi

 \newcommand{\Hom}{\mathrm{Hom}}
 
 \newcommand{\Id}{\mathrm{Id}}

 \newcommand{\cB}{{\mathcal B}}

 \newcommand{\cF}{{\mathcal{F}}}

\newcommand{\cA}{{\mathcal A}}

 \newcommand{\C}{{\mathbb C}}
 
 \newcommand{\R}{{\mathbb R}}
 \newcommand{\Z}{{\mathbb Z}}

 \newcommand{\N}{{\mathbb N}}

 \newcommand{\ve}{{\varepsilon}}

 \def\g{\mathfrak{g}}

 \def\so{\mathfrak{so}}

 \def\a{\alpha}

 \def\e{\epsilon}

 \def\bk{\mathds{k}}

\DeclareFontFamily{OT1}{pzc}{}
\DeclareFontShape{OT1}{pzc}{m}{it}{<-> s * [1.15] pzcmi7t}{}

\newcommand{\Sym}{\mathrm{Sym}}

\newcommand{\UE}{{\mathscr{U}}}

\newcommand{\isomto}{\stackrel{\sim}{\rightarrow}}

\begin{document}

\allowdisplaybreaks

\newcommand{\arXivNumber}{2601.07803}

\renewcommand{\PaperNumber}{092}

\FirstPageHeading

\ShortArticleName{Integration \`a la Harish-Chandra for Bi-Graded Lie Algebras}

\ArticleName{Integration \`a la Harish-Chandra \\ for Bi-Graded Lie Algebras}

\Author{Alexei KOTOV~$^{\rm a}$, Vladimir SALNIKOV~$^{\rm b}$ and Olga SALNIKOVA CHEKERES~$^{\rm bc}$}

\AuthorNameForHeading{A.~Kotov, V.~Salnikov and O.~Salnikova Chekeres}

\Address{$^{\rm a)}$~Faculty of Science, University of Hradec Kr\'alov\'e,\\
\hphantom{$^{\rm a)}$}~Rokitanskeho 62, Hradec Kr\'alov\'e, 50003, Czech Republic}
\EmailD{\mail{oleksii.kotov@uhk.cz}}

\Address{$^{\rm b)}$~LaSIE UMR 7356 CNRS / La Rochelle University, \\
\hphantom{$^{\rm b)}$}~Avenue Michel Cr\'epeau, La Rochelle, 17042, France}
\EmailD{\mail{vladimir.salnikov@univ-lr.fr}}

\Address{$^{\rm c)}$~M\&MoCS Center, Universit\`a degli Studi dell'Aquila, \\
\hphantom{$^{\rm c)}$}~Piazzale Ernesto Pontieri, Monteluco di Roio, L'Aquila, 67100, Italy}
\EmailD{\mail{olya.chekeres@gmail.com}}

\ArticleDates{Received January 27, 2026, in final form September 09, 2026; Published online September 24, 2026}

\Abstract{We study $\mathbb{Z}_2\times\mathbb{Z}_2$ bi-graded Lie algebras. We describe their properties in relation to Lie superalgebras with some compatible structures. Then we focus on the approach to the Lie group-algebra correspondence based on Harish-Chandra pairs and provide some examples of application of it in the bi-graded setting.}


\Classification{17B70; 58A50}

\section{Introduction}

For a \emph{graded algebra}, graded by any abelian group or, more generally, by any abelian monoid, one can say that in the definition of the corresponding category, including its objects and morphisms, there is no conceptual dependence on that group or monoid.

For an abelian monoid $\Gamma$, a \emph{$\Gamma$-graded vector space} over a field $\bk$ is understood as a vector space $V$ with a decomposition
\[
V=\bigoplus_{\gamma\in \Gamma}V_\gamma,
\]
where each $V_\gamma$ is a subspace of homogeneous elements of degree (weight) $\gamma$. Morphisms in the corresponding category are linear maps $\phi\colon V\to W$, either preserving degrees ($\phi (V_\gamma)\subset W_\gamma$), or, more generally, maps of fixed degree $\delta\in \Gamma$ satisfying $\phi (V_\gamma)\subset W_{\gamma+\delta}$.

 Graded vector spaces naturally form a \emph{tensor $($or monoidal$)$ category}, which makes it possible to define an \emph{associative graded algebra} within this category. The presence of a compatible symmetry turning the category of graded vector spaces into a symmetric tensor category further enriches it by introducing the notion of \emph{commutative associative graded algebras}. This symmetry is often determined by a sign rule.

\emph{Graded geometry} extends this framework to manifolds or more general spaces equipped with analogous graded structures. However, it does not possess a single canonical approach, owing to variations in conventions concerning coordinates, admissible functions, and sign rules. The notion of an ``analogous structure'' may admit different interpretations.

 Perhaps the only exception is the case of super structures, which are modeled by $\Z_2$-graded spaces equipped with a symmetry of the form \smash{$v\otimes w \mto (-1)^{p(v)p(w)} w\otimes v$}, where $v\in V$, ${w\in W}$ for $\Z_2$-graded spaces $V$ and $W$, and $p$ denotes the parity of an element, taking values in~${\Gamma \!=\! \Z / 2\Z}$. Here, the global geometric object, a \emph{supermanifold}, is a topological space endowed with a structure sheaf or a \emph{locally ringed space} modeled by a free algebra generated by homogeneous elements. By a free algebra we mean an algebra generated by smooth (analytic, \dots) functions in even (i.e., of degree 0) variables and polynomial expressions in odd, or \emph{Grassmann}, variables (of degree~1). Morphisms between such objects are morphisms of the corresponding locally ringed spaces (cf.~\cite{Leites80,AAVoronov84}).

 We can continue to play the same game, namely, to use locally ringed spaces for other gradings. For the category of vector spaces graded by the natural numbers (including zero), one should then consider the sheaf of functions locally generated by smooth (analytic, algebraic, etc.) functions in variables of degree zero and by polynomials in variables of positive degree. For the case of the group of integers, polynomials in generators of nonzero degree are no longer sufficient, and the polynomial algebra has to be completed with respect to a certain natural filtration. The corresponding category, called \emph{semi-formal manifolds}, is described in \cite{KS2024}, following the idea of a filtration in the affine case from \cite{felder-kazhdan} (an equivalent in meaning definition is used in \cite{Smolka-Vysoky2025, Vysoky2022}).

 Another approach to $\N$-graded manifolds was proposed in \cite{grab-hom}: instead of a ringed space, one considers a manifold equipped with an action of a multiplicative monoid $\R_{\ge 0}$. This point of view on graded manifolds has found interesting applications in various problems of differential geometry. The infinitesimal generator corresponding to the monoid action takes the form of an~Euler vector field. This means that locally on the manifold one can choose smooth coordinates which are homogeneous with respect to this field and have nonnegative integer (or, more generally, nonnegative real) weights. Moreover, such a manifold carries a canonical structure of a~graded bundle over the set of stable points, which forms a smooth submanifold.

 This approach extends to
$\Z$-graded manifolds as follows.\footnote{Historically, it was proposed even earlier in \cite{Voronov2002}; see applications in \cite{GrabowskaGrabowski2024, HCP-instances} and a more detailed review on graded geometry in \cite{voronov19}).} A $\Z$-graded manifold is defined to be a smooth (super)manifold equipped with a vector field, also referred to as the Euler field, which admits homogeneous smooth coordinates with integer (or, more generally, real) weights in a neighborhood of the points of its zero locus. In particular, it follows that the set of fixed points of the Euler field forms a smooth submanifold. The corresponding category is straightforward to define by requiring morphisms to be morphisms of the underlying supermanifolds that map one vector field to another. The semi-formal $\Z$-graded manifolds in the sense of \cite{KS2024} then turn out to be isomorphic to the formal neighborhood of the zero locus (see the examples in \cite{GrabowskaGrabowski2024, HCP-instances}).
Let us emphasize that the grading in the above sense, defined by means of the Euler vector field, does not fix the sign rule for multiplication. In particular, a \(\Z\)-graded manifold can be even, that is, its algebra of smooth functions can be commutative in the usual sense. Another particular case is when the superparity of the coordinates is given by the reduction modulo \(2\) of their integer weights. In general, however, neither of these options is required; we assume only the existence of coordinates with both pure superparity and pure integer weights. The advantage of the approach based on Euler vector fields is that morphisms in the corresponding category can be naturally defined as morphisms of supermanifolds that intertwine the Euler vector fields.

 The case of a \(\Gamma = \mathbb{Z}_2 \times \mathbb{Z}_2\) grading (we will call it \emph{bi-grading} in the rest of the article, if it does not lead to any confusion) is of interest both from a mathematical standpoint and in applications (see Survey \cite{Toppan2022talk} and the references therein, also \cite{CheSa23, Lychagin1995, Genoud-Ovsienko2010, VS-aksz}). Vector spaces endowed with such a grading decompose into components of degrees \((0,0)\), \((1,0)\), \((0,1)\), and \((1,1)\). The grading of each homogeneous element $v$ has the form $\ve (v) = (\ve_1 (v), \ve_2 (v))$. It is straightforward to define the corresponding tensor category. The symmetry in this category can be specified in different ways. There are two most ``popular'' and natural ones:

\emph{Method $1$} (Deligne sign convention): $ v\otimes w \mto (-1)^{\ve (v)\ve (w)} w\otimes v$, where
\begin{equation}\label{intro:Deligne}
 \ve (v)\ve (w):= \ve_1 (v)\ve_1 (w) + \ve_2 (v)\ve_2 (w).
\end{equation}

\emph{Method $2$} (Bernstein--Leites sign convention):
$v\otimes w \mto (-1)^{\ve (v)\ve (w)} w\otimes v$, where
\begin{equation*}
\ve (v)\ve (w):= (\ve_1 (v) + \ve_2 (v))(\ve_1 (w)+\ve_2 (w)).
\end{equation*}

This leads to two distinct definitions of a commutative algebra. Nevertheless, the formula
\begin{equation}\label{intro:twist}
 v\otimes w \mto (-1)^{\ve_1 (v)\ve_2 (w)} v\otimes w
\end{equation}
for all elements $v$ and $w$ establishes an equivalence between the corresponding categories. More precisely, the second category ``forgets'' the double grading when computing the sign, responding only to the total parity of the elements. However, in order to define morphisms correctly, one must require them to commute with both gradings. In this context, one should speak of a symmetric tensor category of super vector spaces equipped with a compatible involutive automorphism.\footnote{The more general case of \(\mathbb{Z}_2^n\)-gradings can be treated in a similar manner. However, the resulting theory is more complicated and, in our opinion, less elegant. In any case, this topic lies beyond the scope of the present paper.}

 In defining \(\mathbb{Z}_2 \times \mathbb{Z}_2\)-graded ``global'' spaces, a ``semi-formal'' approach is traditionally used (see \cite{BrIbPo2020, BrIbPo2021,co-gra-po}). In this construction of the function sheaf, only the coordinates of degree \((0,0)\) are allowed to have smooth functions, while the remaining coordinates correspond to formal power series in the generators of degrees \((1,1)\) or Grassmann polynomials in the generators of degrees~\({(1,0)}\) and \((0,1)\).

 On the other hand, instead of working with ringed spaces, one may consider two alternative categories.\footnote{A.K. is grateful to Elizaveta Vishnyakova for drawing attention to the paper \cite{Santamaria2026} (see also \cite{Vishnyakova2019}), which explains that the sign rule can be extended from formal power series to smooth functions.}

\emph{Category $1$}: Smooth manifolds equipped with an involutive automorphism and a ``sheaf'' of commutative functions in the sense of the sign rule of type 1. The word ``sheaf'' is placed in quotation marks here because all operations can be defined only on open subsets and coverings that are invariant under the involution.

\emph{Category $2$}: Smooth supermanifolds endowed with an involutive automorphism.

 In contrast to the algebraic case, for smooth (super)manifolds the procedure~\eqref{intro:twist} defined above does not yield an equivalence between symmetric tensor categories~1 and~2, since the corresponding operation ``pushes'' one outside the class of smooth functions (see Section~\ref{app:smooth}).

 This paper continues our series of works \cite{Bonavolonta:1304.0394,Guarin-Kotov_2026,DGLG, KLS2023, HCP-instances,KS2024} related to graded geometry in general and, in particular, the problem of integration in the setting. In~\cite{DGLG}, we have solved the Lie group--algebra correspondence problem for differential $\N$-graded Lie groups, in terms of equivalence of categories with differential graded Harish-Chandra pairs. This study was extended to the $\Z$-graded case and some more advanced examples were considered in \cite{HCP-instances} after having properly defined the category of (semi-formal) $\Z$-graded manifolds in \cite{KS2024}.
In that same paper~\cite{KS2024} we have mentioned in passing that most of the constructions make sense for any $\Gamma$-graded manifold, $\Gamma$ being a cancellative monoid. In the current study we make this remark explicit for the first near-at-hand non-trivial example of $\Gamma = \Z_2 \times \Z_2$.

The rest of this paper is organized as follows.

We start Section~\ref{sec:bigraded_Lie} by a brief description of the corresponding category with a couple of natural examples, in view of the above general philosophy. In Section~\ref{sec:super_vs_super-super},
we describe the link between bi-graded and super algebras with extra structure. We start calling the former ones ``super-super'' and explain why this is appropriate.
The shown equivalence of categories is rather straightforward on the objects side, but slightly subtle for the morphisms.
In Section~\ref{sec:HCP}, we present the Harish-Chandra construction of group objects for bi-graded Lie algebras with the sign convention given by the equation \eqref{intro:Deligne} above. We pay attention to the subtleties and features of the multiple grading; in particular we explore the integration of the maximal even subalgebra in the context of Harish-Chandra pairs. Finally, in Section~\ref{app:smooth}, we briefly describe the category of noncommutative spaces (that is, spaces that are commutative according to the~Deligne sign convention) that satisfy the notion of smoothness introduced there and to which these group objects belong. We show that, unlike in the algebraic case, the operation~\eqref{intro:twist} above does not yield an equivalence with the category of smooth supermanifolds equipped with an involution.

\section[(Z\_2 times Z\_2))-graded Lie algebras]{$\boldsymbol{(\Z_2\times \Z_2)}$-graded Lie algebras}\label{sec:bigraded_Lie}

By $(\Z_2\times\Z_2)$-graded vector spaces we mean vector spaces of the form
\[ V = \bigoplus V_{i,j}\] for the indices $i, j \in \Z_2 \equiv \Z/2\Z$.
We denote the grading on these spaces as
\[ \ve(\cdot) = (\ve_1(\cdot), \ve_2(\cdot)) \in \Z_2\times\Z_2.\]

$(\Z_2\times\Z_2)$-graded vector spaces form a tensor category with respect to the standard tensor product:
\begin{gather*}
V\otimes W = \bigoplus_{\gamma\in \Z_2\times \Z_2} \big( V\otimes W\big)_\gamma,
\qquad
\big( V\otimes W\big)_\gamma := \bigoplus_{\gamma_1+\gamma_2=\gamma} V_{\gamma_1}\otimes W_{\gamma_2},
\end{gather*}
where $\gamma, \gamma_1, \gamma_2\in \Z_2\times\Z_2=\Z_2^2$.

For the commutation relations, we adopt the Deligne sign convention for the symmetry $V\otimes W\simeq W\otimes V$ of the tensor (monoidal) category (cf.\ \cite{Etingof2015}) of $(\Z_2\times\Z_2)$-graded vector spaces
\begin{gather}\label{eq:Deligne}
\tau\colon\ v\otimes w \mto (-1)^{\ve(v)\ve(w)} w \otimes v,
\end{gather}
where $v\in V$, $w\in W$, and $\ve(v)\ve(w) \colon= \ve_1(v)\ve_1(w) + \ve_2(v)\ve_2(w)$.

\begin{Remark}
We will call $\ve(\cdot)$ \emph{bi-degree} to stress its nature. We will also use the wording \emph{bi-graded} when it does not lead to confusion.
\end{Remark}

\begin{Definition}[\cite{Rittenberg1978,Rittenberg-Wyler1978,Scheunert1979}]
 A \emph{$(\Z_2\times\Z_2)$-graded Lie algebra} is a~$(\Z_2\times\Z_2)$-graded vector space~$\g$ equipped with a bilinear operation $[\cdot, \cdot] \colon \g \times \g \to \g$ of bi-degree $(0, 0)$, satisfying the following two properties:
 \begin{itemize}\itemsep=0pt
 \item bi-graded antisymmetry:
 $[a, b] = - (-1)^{\ve(a)\ve(b)} [b, a]$;
\item viewed as a derivation parametrized by the first element $[a, \cdot]$, bi-graded Leibniz identity of bi-degree $\ve(a) =(\ve_1(a), \ve_2(a))$:
 \[ [a, [b, c]] = [[a, b], c] + (-1)^{\ve(a)\ve(b)}[b, [a, c]].\]
 \end{itemize}
\end{Definition}

\begin{Remark}
 With the bi-graded anti-symmetry, it is easy to check that the bi-graded Leibniz identity can be rewritten as the familiar graded Jacobi identity:
 \[
 (-1)^{\ve(a)\ve(c)} [a,[b,c]] +
 (-1)^{\ve(c)\ve(b) } [c,[a,b]] +
 (-1)^{\ve(b)\ve(a) } [b,[c,a]] = 0.
 \]
 Note that here, as well as the definition above, the formula graphically resembles the usual super case, but conceptually the definition of $\ve(\cdot) \ve(\cdot)$ is crucial.
\end{Remark}

\begin{Example}
An example of a $(\Z_2\times\Z_2)$-graded Lie algebra comes naturally from the example of bi-graded operators on bi-graded vector spaces mentioned above: it is natural to consider the bracket given by the bi-graded commutator
\[ [a, b] := a\cdot b - (-1)^{\ve_1(a)\ve_1(b) + \ve_2(a)\ve_2(b)} b \cdot a.
\]
The same construction holds more generally for any associative product ``$\cdot$'' in a $(\Z_2\times \Z_2)$-graded associative algebra $C$. From this point onward, we will refer to this example as the (canonical) bi-graded Lie algebra corresponding to \(C\).
\end{Example}

\begin{Remark} \label{rem:g+}
 As above, decompose $\g$ as a bi-graded vector space
 \[ \g \equiv \g_{00} \oplus \g_{01} \oplus \g_{10} \oplus \g_{11},
 \]
then the following properties are direct consequences of the definition:
\begin{itemize}\itemsep=0pt
 \item $\g_+ := \g_{00} \oplus \g_{11}$ can be viewed as an even $\Z_2$-graded Lie algebra.
 \item $\g_{00} \oplus \g_{01}$ and $\g_{00} \oplus \g_{10}$ are Lie super algebras in the usual meaning of the term.
 \item {The brackets between $\g_{01}$ and $\g_{11}$, as well as those between $\g_{10}$ and $\g_{11}$, are symmetric, whereas the brackets between $\g_{01}$ and $\g_{10}$ are skew-symmetric.}
\end{itemize}
\end{Remark}

\begin{Example}[unitary structure]
 Let $A=A_0\oplus A_1$ be an associative complex $\Z_2$-graded algebra with a compatible unitary structure, that is, an anti-linear involution $*$, which commutes with the grading automorphism. More precisely, $*^2=\Id$,
$*(A_0)\subset A_0$, $*(A_1)\subset A_1$,
and
\begin{gather*}
*(\a a) = \bar\a *(a)
\qquad\forall \a\in\C,\ a\in A, \qquad
*(ab) = *(b)*(a)
\qquad \forall a,b\in A.
\end{gather*}
Let
\begin{gather*}
u_r = \{ a\in A_r \mid *(a)=-a\}, \qquad
h_r = \{ a\in A_r \mid *(a)=a\}
\end{gather*}
for $r=0,1$.

We equip the subspaces $u_r$, $h_r$, $r=0,1$ with the following bi-grading:
\[
\ve (u_r)=(r,0), \qquad \ve (h_r)=(r+1 \ (\bmod~2),1).
\]
The structure of a bi-graded Lie algebra on $A$ is given by a bilinear operation defined by the formulas
\[
[a,b] =\begin{cases}
[a,b]_-, & a,b\in A_0 \ \mathrm{or}\ a\in u_1, b\in h_1,\\
{\rm i}[a,b]_+, & a\in u_1,\ b\in h_0 \ \mathrm{or}\ b\in u_1,\\
-{\rm i}[a,b]_+, & a\in h_1,\ b\in h_0 \ \mathrm{or}\ b\in h_1,
\end{cases}
\]
where $[a,b]_-=ab-ba$ and $[a,b]_+=ab+ba$. It is assumed everywhere that the operation is skew-symmetric with respect to the bi-grading defined above.

 It is obvious that the operation defined in this way respects the bi-grading. As an exercise, the reader can directly verify the bi-graded Jacobi identity through computation. We will prove this last property in the following example by embedding the resulting algebra as a subalgebra into another bi-graded Lie algebra.
\end{Example}

\begin{Example}
Consider the following commutative (in the usual sense) associative algebra:
$B = \{x_0 + x_1 q_1 + x_2 q_2 + x_3 q_3, $
 with $x_i \in \C$ and
 with the relations $q_1^2={\rm i}$, $q_2^2=-{\rm i}$, $q_1 q_3={\rm i} q_2$, $q_2 q_3=-{\rm i} q_1 \}$. This automatically implies
 $q_3^2=1$ and $q_1q_2=q_3$. Let us assign the following bi-gradings to the generators of the algebra: $\ve (1)=(0,0)$, $\ve(q_1)=(1,0)$, $\ve (q_2)=(0,1)$, and $\ve (q_3)=(1,1)$. This endows the algebra with a bi-grading compatible with the multiplication.

The algebra $B$ is isomorphic\footnote{We thank an anonymous referee for their suggestion of this simplification.} to $\C\otimes_{\R}\C$ via the correspondence
\[
q_1=\lambda^3({\rm i}\otimes 1), \qquad
q_2=\lambda(1\otimes {\rm i}), \qquad
q_3=-({\rm i}\otimes {\rm i}), \qquad\text{where } \lambda=\exp\biggl(\frac{\pi {\rm i}}{4}\biggr).
\]

 We associate to any $\Z_2$-graded associative algebra $C=C_0\oplus C_1$ a $(\Z_2\times \Z_2)$-graded associative algebra
$\widetilde{ C}=C_0 \oplus C_0 q_3\oplus C_1 q_1 \oplus C_1 q_2$.
The bi-graded Lie algebra from the previous example admits the following embedding into the bi-graded Lie algebra corresponding to $\widetilde{C}$:
\[
u_0 \oplus u_1 q_1 \oplus h_1 q_2 \oplus h_0 q_3.
\]
\end{Example}

\section{Super-super vs super}\label{sec:super_vs_super-super}

 Note that every $(\Z_2\times \Z_2)$-bi-graded vector space is $\Z_2$-graded with respect to the total grading $p(\cdot):= \ve_1 (\cdot)+\ve_2 (\cdot)$, which we will also refer to as \emph{parity}. It is clear that the parity of an element contains only half of the information about its bi-grading. The other half can be recovered by means of an involutive automorphism $\sigma$ defined as follows:
 \begin{gather}\label{eq:automorhism}
 \sigma (v):=(-1)^{\ve_2 (v)}v.
 \end{gather}
 It is now easy to see that the category of bi-graded vector spaces is equivalent to the category of vector superspaces equipped with an involutive automorphism. Objects in the latter category are vector superspaces equipped with an involutive automorphism, and morphisms are, accordingly, the equivariant morphisms between such vector superspaces.

 On the other hand, the symmetry in the tensor category of superspaces is given by the isomorphism
$
v\otimes w \mto (-1)^{p(v)p(w)} w \otimes v$,
where $v\in V$, $w\in W$. In order for the equivalence of categories defined above to become an~equivalence of the corresponding symmetric tensor categories as well, one needs to supplement the functor defining the equivalence with an isomorphism
\begin{gather}\label{eq:unbraiding}
v\otimes w \mto (-1)^{\ve_1 (v)\ve_2 (w)} v\otimes w,
\end{gather}
where $v$, $w$ are as above. By applying the technique defined above to a bi-graded Lie algebra~$\g$, we obtain a Lie superalgebra whose underlying space is $\g$, whose ``new'' superbracket, as one can easily guess, is given by the formula
\begin{gather}\label{eq:unbrading_bracket}
[a, b]_{s}:= (-1)^{\ve_1 (a)\ve_2 (b)} [a,b]
\end{gather}
for all $a,b\in\g$ of pure bi-degree.

\begin{Proposition}
 We obtain an equivalence between the category of bi-graded Lie algebras and the category of Lie superalgebras equipped with an involutive\footnote{Similar involutions have been studied in theliterature, including \cite{chuah, fioresi, pellegrini, serganova}.} automorphism. Morphisms in the latter category are those morphisms of the corresponding Lie superalgebras that interchange the involutions.
\end{Proposition}

\begin{proof} As follows from general theory, the resulting operation defines a structure of a Lie superalgebra. One can verify explicitly that
\[
(-1)^{\ve(a_1)\ve(a_3)} [a_1,[a_2,a_3]] + {\rm c.p.}
 =(-1)^\a \big( (-1)^{p(a_1)p(a_3)} [a_1,[a_2,a_3]_s]_s + {\rm c.p.}\big),
\]
where $\a = \ve_1 (a_1)\ve_2 (a_2)+\ve_1 (a_2)\ve_2 (a_3)+\ve_1 (a_3)\ve_2 (a_1) $. This means that one Jacobi identity holds if and only if the other one does.

 Moreover, the correspondence is functorial in the sense that morphisms of the algebra $\g$ are in one-to-one correspondence with morphisms of the resulting Lie superalgebra that preserve the additional grading. Indeed, let $\phi\colon \g_1\to \g_2$ be a bi-graded Lie algebra morphism, that is, a~bi-degree preserving map which satisfies $\phi [a,b]=[\phi(a),\phi (b)]$ for all $a,b\in\g_1$. Then
\[
\phi [a,b]_s=(-1)^{\ve_1 (a)\ve_2 (b)} \phi [a,b] =
(-1)^{\ve_1 (a)\ve_2 (b)} [\phi (a),\phi (b)]=
[\phi(a),\phi (b)]_s.
\]
The converse statement holds as well: any morphism of the corresponding Lie superalgebras that respects the bi-grading canonically determines a morphism of bi-graded Lie algebras. Moreover, under this correspondence, the composition of morphisms on either side naturally corresponds to the composition on the other. In particular, the involutive automorphism \eqref{eq:automorhism} of the bi-graded algebra $\g$ induces an involutive automorphism of the corresponding Lie superalgebra. It is obvious that a Lie superalgebra morphism preserves the additional grading if and only if it is equivariant with respect to the automorphisms.
\end{proof}

\begin{Remark}
 Suppose that the bi-graded Lie algebra contains only the components $(0,0)$ and~$(1,1)$; that is, it is an ``ordinary'' (even) Lie algebra endowed with $\Z_2$-grading, or equivalently, with an involutive automorphism. The decomposition $\g=\g_+=\g_{00}\oplus \g_{11}$
 is referred to in the literature as \emph{symmetric} (such Lie algebras arise in the geometry of \emph{symmetric spaces}, cf.~\cite{Helgason2024}), and its homogeneous components $\g_{00}$ and $\g_{11}$ are called a \emph{Cartan pair}. In this case, transformation \eqref{eq:unbrading_bracket} induces an involutive automorphism of the category of $\Z_2$-graded Lie algebras. It is easy to see that all operations are preserved, except for the bracket on~${\g_{11}\otimes \g_{11}\to \g_{00}}$, which changes sign.
\end{Remark}

\begin{Example}\label{ex:orthogonal}
 Consider the orthogonal Lie algebra $\g=\so (3)$ with the standard basis $\{e_i\}$, $i=1,2,3$, satisfying the relations
 \[
 [e_1,e_2]=e_3, \qquad [e_2,e_3]=e_1, \qquad [e_3,e_1]=e_2
 \]
 The invertible linear $\sigma$ operator defined by formula
 \[
 \sigma (e_1)=e_1, \qquad \sigma (e_2)=-e_2, \qquad \sigma (e_3)=-e_3
 \]
 induces an involutive automorphism of $\g$. It is inner, that is, determined by the adjoint action of an element of the orthogonal group which squares to the identity. Indeed, this is evident for the standard linear representation:
 \begin{gather*}
 e_1 =
 \left(\begin{smallmatrix} 0 & 0 & \hphantom{-}0 \\
 0 & 0 & -1 \\
 0 & 1 & \hphantom{-}0
 \end{smallmatrix}\right)
 , \qquad
 e_2 =
 \left(\begin{smallmatrix} \hphantom{-}0 & 0 & 1 \\
 \hphantom{-}0 & 0 & 0 \\
 -1 & 0 & 0
\end{smallmatrix}\right), \qquad
 e_3 =
 \left(\begin{smallmatrix} 0 & -1 & 0 \\
 1 & \hphantom{-}0 & 0 \\
 0 & \hphantom{-}0 & 0
\end{smallmatrix}\right), \qquad
 \sigma =
 \left(\begin{smallmatrix} 1 & 0 & \hphantom{-}0 \\
 0 & 0 & -1 \\
 0 & 1 & \hphantom{-}0
\end{smallmatrix}\right).
 \end{gather*}
One readily checks that transformation \eqref{eq:unbrading_bracket} sends $\g=\so (3)$ to the Lie algebra $\so (1,2)$.

 Furthermore, an analogous construction applies to any Lie supergroup whose even part contains ${\rm SO}(3)$ as a subgroup (for example, the orthosymplectic group ${\rm OSp}(3,2n)$). The adjoint action of $\sigma$ yields an involutive automorphism of the corresponding Lie algebra and, consequently, an additional $\Z_2$-grading compatible with the superstructure. Transformation \eqref{eq:unbrading_bracket} then sends the associated Lie superalgebra to a bi-graded Lie algebra such that $\so(3)$
 is mapped to $\so(1,2)$ (and vice versa).
\end{Example}

In Section \ref{sec:HCP}, we will show that the group objects integrating a bi-graded Lie algebra and the corresponding Lie superalgebra with an involutive automorphism may have different topological properties. We will see that, at the level of smooth graded and bi-graded spaces, the ``sign trick''~\eqref{eq:unbraiding} does not work quite as literally as one might expect.

\section{Integration and Harish-Chandra pairs} \label{sec:HCP}
Through this whole section, we will again use the decomposition
$\g = \g_{00} \oplus \g_{01} \oplus \g_{10} \oplus \g_{11}$,
and a short-hand notation $\ve(a)\ve(b) \equiv \ve_1(a)\ve_1(b) + \ve_2(a)\ve_2(b)$.

 The formal integration of a Lie algebra is given by its universal enveloping algebra;
i.e., as a~Hopf algebra it determines the resulting formal group object.
For a bi-graded Lie algebra, its construction is rather straightforward: we shall follow the standard procedure.

\begin{Definition}
The \emph{universal enveloping algebra} of $\g$ is the unital associative algebra $\UE(\g)$, defined as the quotient algebra $T(\g)/I$, where $T(\g)$ is the tensor algebra of $\g$ and $I$ is the two-sided ideal generated by all elements of the form
\[a\otimes b-(-1)^{\ve(a)\ve(b)}b\otimes a-[a,b] \qquad\text{for } a,b\in \g.\]
\end{Definition}

As in the classical case, there is a canonical map $\iota\colon \g\to \UE(\g)$. Later on we will see that it is an inclusion, which will allow us to identify $\g$ with its image in $\UE(\g)$ under the above map.

\begin{Remark}
 By the standard arguments, $\UE(\g)$ can be endowed with a Hopf algebra structure as follows.

In addition to the (generally non-commutative) multiplication, $\UE(\g)$ is equipped with
\begin{itemize}\itemsep=0pt
 \item an even compatible bi-degree $(0,0)$ coassociative super cocommutative comultiplication $\Delta\colon \UE(\g)\to \UE(\g)\otimes \UE(\g)$, such that
$\Delta (ab)=\Delta (a)\Delta (b)$
for all $a,b\in \UE(\g)$ and a counit $\UE(\g)\to\bk$, which descend from the ones on $T(\g)$. In~particular, an element $a$ is \emph{primitive}, i.e., it satisfies
$\Delta (a)=a\otimes 1+ 1\otimes a$ if and only if it belongs to $\iota(\g)$.
 \item an antipode $S\colon\UE(\g)\to\UE(\g)$, defined such that
 \[
 S(x_1\dots x_m)=(-1)^m
 \epsilon(x_1,\dots, x_m) x_m\dots x_1 \qquad \forall x_1,\dots,x_m\in\g,
 \]
 where the prefactor \smash{$\epsilon(x_1,\dots, x_m) = (-1)^{\sum_{i<j}\ve(x_i)\ve(x_j)}$}.
\end{itemize}
\end{Remark}

Note that the universal enveloping algebra is cocommutative in the sense of symmetry given by the equation \eqref{intro:Deligne} from the introduction.
For the zero bracket, the corresponding universal enveloping algebra will be the \emph{symmetric algebra} generated by $\g$, which we shall denote by~$\Sym (\g)$. As a special case of a more general construction, the symmetric algebra is also a~Hopf algebra. The symmetric algebra is not only cocommutative but also commutative in the sense of the same symmetry~\eqref{intro:Deligne}.

\begin{Proposition}[Poincar\'{e}--Birkhoff--Witt theorem]
Given a basis $\mathcal{B}$ of $\g$ one can construct explicitly a basis of $\UE(\g)$, by a proper ordering of canonical monomials of $\mathcal{B}$.
In particular, the linear map $\iota$ is an inclusion and extends to an isomorphism of graded vector spaces ${\Sym (\g)\isomto \UE(\g)}$. Moreover, by the Weyl formula, the linear isomorphism can be made canonical, so that it also becomes an isomorphism of cocommutative coalgebras:
\[
\Sym (\g)\ni x_1 x_2 \dots x_m \mto \frac{1}{m!} \sum_{\sigma\in S_m}
\e (\sigma) x_{\sigma (1)} x_{\sigma (2)} \dots x_{\sigma (m)} \in \UE (\g)
\]
for all $x_1, x_2, \dots, x_m \in \g$, where $\e (\sigma)$ is the Koszul sign of the permutation $\sigma\in S_m$, generated by the rule \eqref{eq:Deligne}.
\end{Proposition}

\begin{proof}
The details of the proof are no different from those for the Lie superalgebra or the graded Lie algebra (see \cite{HCP-instances}, also \cite{FHT,example}).
\end{proof}

We now follow the strategy from \cite{HCP-instances} to produce an integrating object for $\g$.

Consider the ``even'' subalgebra of $\g$,
$ \g_+ = \g_{00} \oplus \g_{11}$,
which is a Lie algebra in the usual sense (see Remark~\ref{rem:g+}).
Ignoring the multiple gradings, it can be integrated to a Lie group $G_+$ equipped with an involutive automorphism.
The following proposition holds.

\begin{Proposition}\label{prop:G+_integration}
Every Lie algebra admitting an involutive automorphism integrates to a Lie group equipped with a compatible involutive automorphism. The integration is unique when the Lie group is required to be simply connected. Moreover, the adjoint action of the Lie subalgebra~$\g_+$ integrates to a representation of this Lie group that is compatible with both the bi-graded Lie bracket and the involutive automorphism.
\end{Proposition}

\begin{proof} The proof relies on the standard integration technique for Lie algebras and their representations. First, construct the group and representation locally near the identity using the Baker--Campbell--Hausdorff formula for the product. Extend the result analytically along paths. For the simply connected Lie group, the extension is unique; hence, compatibility with both the bi-graded Lie bracket and the involutive automorphism follows automatically.
\end{proof}

Consider now the pair $(\g,\g_+)$ as a (generalized) Harish-Chandra pair within the Harish-Chandra approach \cite{Kostant:1975, Vishnyakova:2011}; that is, replace it with the pair $(G_+,\g)$ as in Proposition~\ref{prop:G+_integration}.

\begin{Proposition}\label{prop:bigraded_HC}
Consider the space
$\cF(G):=
\Hom_{\UE(\g_+)}(
\UE(\g),\cF(G_+)
)$,
where $\cF(G_+)$ is the algebra of functions on $G_+$. The following statements hold true:
\begin{enumerate}\itemsep=0pt
 \item[$(1)$] $\cF(G)$ is a bi-graded Hopf algebra $($commutative in the sense of the Deligne sign convention~$\eqref{intro:Deligne})$, which we regard as the algebra of functions on~$G$.

 \item[$(2)$] As a bi-graded vector space, $\cF(G)$ is isomorphic to the space of sections of the trivial vector bundle
 $ G_+\times \Lambda^{\bullet} (\g_{10}^*\oplus\g_{01}^*)$.

 \item[$(3)$] The integration procedure in {\rm \cite{example}}, based on the Harish-Chandra pair $(G_{00},\g)$, where $G_{00}$ is the Lie group integrating~$\g_{00}$, yields a structure isomorphic to the formal neighbourhood of~$G_{00}$.
\end{enumerate}
\end{Proposition}

\begin{proof}
The proof does not differ substantially from the standard one for Harish-Chandra pairs and follows the general idea of \cite{Kostant:1975}.
 In particular, commutativity follows from the following simple observation: given a cocommutative coalgebra $\cA$
 and a commutative algebra $\cB$ (in the sense of the symmetry \eqref{intro:Deligne}), the space of linear maps $\cA\to \cB$
 is a commutative algebra with the multiplication
 \[
 (\phi*\psi )(a):=\mu_{\cB}\circ (\phi\otimes \psi )\circ \Delta_{\cA} (a)
 \]
 for each $a\in \cA$, where $\Delta_{\cA}$ and $\mu_{\cB}$ are the comultiplication and multiplication in $\cA$ and $\cB$, respectively. From $\tau\circ \Delta_{\cA}=\Delta_{\cA}$ (cocommutativity of $\cA$) and $\mu_{\cB}\circ\tau =\mu_{\cB}$ (commutativity of~$\cB$), immediately obtain that the resulting algebra is commutative. The action of the universal enveloping algebra of the subalgebra $\g_+$ is compatible with the bigrading. Thus, the invariant subspace inherits the required property of being a bigraded commutative algebra.

 The last statement is proved analogously to \cite[Proposition~2]{HCP-instances}.
\end{proof}

\begin{Remark}
The pair $(G_+, \cF(G))$ from the above Proposition~\ref{prop:bigraded_HC}
represents a Lie group in the specific category of bi-graded manifolds integrating~$\g$, that we briefly discuss in the beginning of Section~\ref{app:smooth}. We address some features of this category and elaborate on the general construction in a separate paper.
\end{Remark}

\begin{Remark}
Despite the naturalness of the Harish-Chandra approach, there is an important point that we would like to emphasize once again: the choice of $\g_+$ makes the resulting construction essentially different from the Harish-Chandra construction of \cite{example}, which is based on the pair $(G_{00},\g)$. The latter yields a group object in the category studied, for example, in \cite{BrIbPo2020,BrIbPo2021, co-gra-po}.

In Proposition~\ref{prop:bigraded_HC}, we follow the general philosophy of \cite{HCP-instances}, namely, to integrate the largest possible ``even'' Lie subalgebra. In the present setting, this subalgebra is $\g_+$. The resulting group object belongs to a different category, which we briefly describe in Section~\ref{app:smooth} (see also Remark~\ref{rem:bigraded_group} below).
\end{Remark}

\begin{Remark}\label{rem:bigraded_group} \
\begin{itemize}\itemsep=0pt
\item The automorphism $\sigma$ acts on $\cF(G)$ compatibly. Moreover, the bi-grading corresponds bijectively to the pair consisting of the first grading and the $\sigma$-action.
 \item Local functions on $G$ are a sheaf over $G_+$
 in the sense that functions restrict to any $\sigma$-invariant subset. The gluing property holds for $\sigma$-invariant (also called \emph{good}) open coverings.
 \item Since Proposition~\ref{prop:bigraded_HC} provides a general construction for any finite-dimensional $(\mathbb Z_2 \times\mathbb Z_2)$-graded Lie algebra,
 examples of such groups can be easily obtained from such Lie algebras, which in turn are in one-to-one correspondence (cf.\ \eqref{eq:unbrading_bracket}) with super Lie algebras equipped with an involution.
 \end{itemize}
\end{Remark}

\section{Some words about super-super manifolds}\label{app:smooth}

In the previous Section~\ref{sec:HCP}, through integration, we obtained a group object of the following form. Suppose $M$ is a smooth manifold with an involution~$\sigma$, and $E\to M$ is an equivariant bundle over $M$
in the sense that $\sigma$ lifts to an involutive bundle isomorphism $\tilde{\sigma}$ covering $\sigma$.
The space of sections of the vector bundle $E^*$ is canonically bigraded: the first grading is induced by the involution $\tilde{\sigma}$, while the second is induced by the canonical involution of the vector bundle, acting by multiplication by $-1$ on each fiber (see the discussion at the beginning of Section~\ref{sec:super_vs_super-super}). The space $\Gamma (E^*)$ generates a bigraded commutative algebra, which we view as the algebra of functions on a globally bigraded manifold with base $M$. More generally, a bigraded manifold is locally modelled by the above construction. Such bigraded manifolds constitute a category, and the integration of a bigraded Lie algebra yields a group-like object in this category.

The next example illustrates the following phenomenon: the equivalence \eqref{eq:unbraiding} between bi-graded and graded-with-involutive-automorphism categories fails for smooth manifolds. More precisely, applying formula \eqref{eq:unbraiding} to the algebra of smooth functions on a manifold equipped with an involutive automorphism yields an algebra belonging to a different category.
\begin{Example}\label{ex:real_line}
 Consider the real line $\R$ equipped with an involution $x\mto -x$, where $x$ is the affine coordinate on $\R$, which induces an involutive automorphism $\sigma$ of the algebra of smooth functions $C^\infty (\R)$. This algebra thus decomposes into the direct sum of $\sigma$-eigenspaces:
 \[
 C^\infty (\R) =C^\infty (\R)_+\oplus C^\infty (\R)_-,
 \qquad C^\infty (\R)_\pm =\{ f\in C^\infty (\R) \mid f(-x)=\pm f(x) \}.
 \]
 Formula \eqref{eq:unbraiding} deforms the pointwise multiplication of smooth functions on the real line as follows:
 \[
 \big(f_++f_-\big)*\big( h_+ +h_-\big)=f_+h_+ -f_-h_- + f_+h_- +f_-h_+,
 \]
where $f_\pm, h_\pm\in C^\infty (\R)_\pm$.
Identifying $C^\infty (\R)$ with
\[
\cA=\{g=f_++ {\rm i}f_- \mid f_\pm \in C^\infty (\R)_\pm\}
\]
as vector spaces, the deformed multiplication yields an algebra isomorphic to the algebra of smooth complex-valued functions on the real line that satisfies
$g(-x)=\overline{g(x)}$.

Clearly, the resulting algebras are not isomorphic: while the spectrum of $C^\infty (\R)$ with finite-dimensional residue fields (always isomorphic to $\R$) is parametrized by $R$, the spectrum of the second algebra $\cA$ is parametrized by $a\in\R_{\ge 0}$ with residue fields $\R$ at $a=0$ and $\C$ at $a>0$.
\end{Example}

Example \ref{ex:real_line} admits a natural and straightforward generalization to the case of an arbitrary smooth manifold $M$ with a (nontrivial) involution $\sigma$. The algebra of smooth real functions on this manifold, equipped with the multiplication deformed according to formula \eqref{eq:unbraiding}, is not isomorphic to the algebra of smooth real functions -- neither on this manifold nor on any other manifold. The latter statement is verified by the same trick involving the representation of functions in the form $f_++{\rm i}f_-$, where $\sigma (f_\pm)=\pm f$ (we identify $\sigma$ with its pullback map on functions), and by counting the finite-dimensional residue fields with respect to maximal ideals: unlike the algebra of real functions, for which the residue field is always $\R$, in the case of the deformed multiplication it is $\R$ at the points stable under the involution and $\C$ at the generic~points.

This indirectly explains why, in Example~\ref{ex:orthogonal}, integrating the Lie algebra before and after applying operation~\eqref{eq:unbrading_bracket} ($\so(3)$ and $\so (1,2)$, respectively) yields topologically distinct Lie groups: in the first case, a compact group; in the second, a non-compact one. Indeed, if operation \eqref{eq:unbraiding}, applied to the algebra of smooth functions on a manifold with an involution, preserved the class of smooth functions, then the group object corresponding to the Lie algebras obtained from each other via \eqref{eq:unbrading_bracket} would be expected to yield function algebras related by the operation \eqref{eq:unbraiding}. This hypothetical property obviously fails, rendering the global category of bi-graded manifolds a richer and more interesting subject of study.

\subsection*{Acknowledgments}
We are thankful to Elizaveta Vishnyakova, Katarzyna Grabowska, Janusz Grabowski for inspiring discussions and comments at various stages of this work. In particular, A.K. is grateful to J. Grabowski for his attention and insightful questions during the visit to IMPAN (Warsaw) in November 2025, without which a common understanding of the category of global bi-graded manifolds would not have been reached; a paper on this more general topic will be published later.
We thank the anonymous referees for their valuable suggestions, relevant questions and remarks that permitted to improve the presentation.

The research of A.K.\ was supported by the grant ``Graded differential geometry with applications'' GA\u{C}R 24-10031K of the Czech Science Foundation. A research visit of O.C.\ and V.S.\ to Hradec Kr\'alov\'e in Spring 2024, partially supported by the same grant GA\u{C}R 24-10031K, provided the setting in which the idea of pursuing the present line of research first emerged. Another meeting related to this work, in Warsaw in 2025, was funded by the WAVE-UNISONO grant.
A research visit of V.S.\ in Summer 2025 was supported by ACI ``Mobilit\'e Internationale'' of the La Rochelle University.

\pdfbookmark[1]{References}{ref}
\LastPageEnding

\end{document}